\documentclass{article}

\usepackage{amsmath,amssymb,amsthm}

\theoremstyle{plain}\newtheorem{theorem}{Theorem}[section]
   \newtheorem{lemma}[theorem]{Lemma}
   \newtheorem{proposition}[theorem]{Proposition}
   
\theoremstyle{definition}\newtheorem{definition}[theorem]{Definition}
\theoremstyle{remark}\newtheorem{remark}[theorem]{Remark}

\begin{document}

%
%
%
%
\def\R{\mathbb{R}}
%
%
\def\N{\mathcal{N}}
\def\M{\mathcal{M}}
%
%
\def\ssm{\smallsetminus}
\def\symdif{\triangle}
\def\Pow#1{\mathcal{P}(#1)}
\def\restrictedto%
{\mathchoice%
    {\!\upharpoonright\!}%
    {\!\upharpoonright\!}%
    {\upharpoonright}%
    {\upharpoonright}%
}
\def\isomorphic{\simeq}
\def\cf{\operatorname{cf}}
\def\op#1#2{\langle #1,#2\rangle}
\def\dom{\operatorname{dom}}
\def\ran{\operatorname{ran}}
\def\st{:}
\def\card#1{\lvert#1\rvert}
\def\size#1{\lvert#1\rvert}
\def\lh#1{\lvert#1\rvert}
\def\Seq#1{\left\langle{#1}\right\rangle}
\def\rk#1{\operatorname{rank}(#1)}  
\def\domn{\leq^*}
%
%
\newcommand{\forces}[1][{}]{\Vdash_{#1}}
\def\V{\mathbf{V}}
\def\tsum{\textstyle\sum}
\def\tbigcup{\textstyle\bigcup}
\def\tbigcap{\textstyle\bigcap}

\def\bP{\mathbb{P}}
\def\bQ{\mathbb{Q}}

\def\bN{\mathbb{N}}

\def\bLOC{\mathbb{LOC}}
\def\LOC{\mathbb{LOC}}

\def\cS{\mathcal{S}}
\def\cT{\mathcal{T}}

\def\Gdelta{$G_\delta$}
\def\Fsigma{$F_\sigma$}

\title{Hechler's theorem for the null ideal}
\author{Maxim R.\ Burke\thanks{%
    Research supported by NSERC. The
    author thanks F.D.~Tall and the Department of Mathematics
    at the University of Toronto for their hospitality during
    the academic year 2003/2004 when the present paper
    was completed.}\\
{\footnotesize Department of Mathematics and Statistics} \\
{\footnotesize University of Prince Edward Island} \\
{\footnotesize Charlottetown PE, Canada C1A 4P3} \\
{\footnotesize burke@upei.ca}
\and
    Masaru Kada\thanks{%
    The author was supported by
    Grant-in-Aid for Young Scientists (B) 14740058, MEXT.}\\
{\footnotesize Department of Computer Sciences}\\
{\footnotesize Kitami Institute of Technology}\\
{\footnotesize Kitami, Hokkaido 090-8507 Japan}\\
{\footnotesize kada@math.cs.kitami-it.ac.jp}
%
}                     
%
%
\maketitle

\renewcommand{\thefootnote}{}
\renewcommand{\thefootnote}{\fnsymbol{footnote}}

\begin{abstract}
We prove the following theorem: For a partially ordered set $Q$
such that every countable subset of $Q$ has a strict upper bound,
there is a forcing notion satisfying the countable chain condition
such that, in the forcing extension, there is a basis of the null
ideal of the real line which is order-isomorphic to $Q$ with
respect to set-inclusion. This is a variation of Hechler's
classical result in the theory of forcing. The corresponding
theorem for the meager ideal was established by Bartoszy\'nski and
Kada.
\end{abstract}

\section{Introduction}\label{sec:intro}

For $f,g\in\omega^\omega$, we say $f\domn g$ if $f(n)\leq g(n)$
for all but finitely many $n<\omega$. The following theorem, which
is due to Hechler \cite{He:cof}, is a classical result in the
theory of forcing (See also \cite{Bur:hechler}).
\begin{theorem}
Suppose that $(Q,\leq)$ is a partially ordered set such that every
countable subset of $Q$ has a strict upper bound in $Q$, that is,
for any countable set $A\subseteq Q$ there is $b\in Q$ such that
$a<b$ for all $a\in A$. Then there is a forcing notion $\bP$
satisfying the countable chain condition such that, in any forcing
extension by $\bP$, $(\omega^\omega,\domn)$ contains a cofinal
subset $\{f_a\st a\in Q\}$ which is order-isomorphic to $Q$, that
is,
\begin{enumerate}
\item for every $g\in\omega^\omega$
    there is $a\in Q$ such that $g\domn f_a$, and
\item for $a,b\in Q$,
    $f_a\domn f_b$ if and only if $a\leq b$.
\end{enumerate}
\end{theorem}
Fuchino and Soukup \cite{FS:chubu2001,So:specnote} introduced the
notion of \emph{spectra}. For a partially ordered set $P$,
\emph{the unbounded set spectrum of $P$} is the set of cardinals
$\kappa$ such that there is an unbounded set in $P$ of size
$\kappa$ without unbounded subsets of size less than $\kappa$.
They also defined several variants of spectra, and investigated
how to manipulate those spectra of $(\omega^\omega,\domn)$ using
Hechler's result. In this context, Soukup asked if the statement
of Hechler's theorem holds for the meager ideal or the null ideal
of the real line with respect to set-inclusion.
Bartoszy\'nski and Kada \cite{BaK:embedm} answered the question
positively for the meager ideal. In the present paper, we will
give a positive answer for the null ideal. These proofs all follow
the same general scheme, but there are substantial technical
difficulties to overcome in our present context resulting in a
more complicated proof.

As Hechler pointed out in \cite{He:cof}, if $Q$ is well-founded,
the conclusion of his theorem can be strenghened to say that
whenever $a<b$ in $Q$, not only does $f_b$ dominate $f_a$, it also
dominates all reals constructible from $f_a$ and the set of its
predecessors in the cofinal family, i.e., $g\domn f_b$ for all
$g\in L[\{f_x: x\leq a\}]$. Hjorth has answered a question of
Hechler by showing that this stronger conclusion cannot hold if
$Q$ is not well-founded.
\begin{quote}
(Hjorth) There is no sequence of reals $f_n\in\omega^\omega$ such
that for each $n<\omega$, $g\domn f_n$ for each
$g\in\omega^\omega\cap L[\{f_i\st i>n\}]$.
\end{quote}
(See \cite[Theorem~0.5 and preceding discussion]{Br:mad}.) These
results carry over to our context. When $Q$ is well-founded, the
proof of Hechler's theorem for the null ideal $\N$ of $[0,1]$
provides a cofinal family $\{H_a\st a\in Q\}$ of Borel sets in
$\N$ so that $(\{H_a\st a\in Q\},\subseteq)$ is isomorphic to $Q$
and, moreover, if $a<b$ then for every Borel null set $A$ coded in
$L[\{H_x: x\leq a\}]$, we have $A\subseteq H_b$. Moreover,
Hjorth's result implies that this stronger conclusion cannot hold
if $Q$ is not well-founded. This can be seen as follows. By
\cite[Theorems 2.2.2 and 2.3.1]{BaJ:set} and their proofs (which
show that the Tukey maps they provide are definable), if
$M\subseteq N$ are transitive models of enough of ZFC and in $N$
there is a Borel null set $B$ such that $A\subseteq B$ for every
Borel null set $A$ coded in $M$, then in $N$ there is an
$f\in\omega^\omega$ which dominates $\omega^\omega\cap M$. Thus,
the existence of a sequence $A_n$ of null Borel sets such that for
each $n<\omega$, $A_n$ contains all null Borel sets coded in
$L[\{A_i\st i>n\}]$ would yield a counterexample to Hjorth's
result.

We will work with the null ideal of the Cantor set $2^\omega$
rather than that of $[0,1]$ or the real line. The distinction
between these spaces is unimportant in our work because there are
Borel isomorphisms between them which preserve null sets.

\section{Combinatorial view of null sets}\label{sec:combnull}

In this section, we review the relationship between Borel null
sets of the Cantor set $2^\omega$ with the standard product
measure and combinatorics on natural numbers, which is described
in \cite{Ba:invmc}.

Choose a strictly increasing function $h\in\omega^\omega$
satisfying $2^{h(n)-h(n-1)}\geq n+1$ for $1\leq n<\omega$ (for
example, just let $h(n)=n^2$). For each $n<\omega$, let
$\{C_i^n\st i<\omega\}$ be a list of all clopen subsets of
$2^\omega$ of measure $2^{-h(n)}$. We assume that such $h$ and
$C^n_i$'s are fixed throughout this paper.

For a function $f\in\omega^\omega$, we define
\[
H_f=\bigcap_{N}\bigcup_{n>N}C^n_{f(n)}.
\]
Then $H_f$ is a \Gdelta{} null set, and every null set $X$ is
covered by $H_f$ for some $f\in\omega^\omega$.

Let $\cS=\prod_{n<\omega}[\omega]^{\leq n}$. We call each
$\varphi\in\cS$ \emph{a slalom}. As in the case of a function, for
a slalom $\varphi\in\cS$ we define
\[
H_\varphi=\bigcap_{N}\bigcup_{n>N}\bigcup_{i\in\varphi(n)}C^n_i.
\]
Then $H_\varphi$ is a \Gdelta{} null set, and the following hold:
\begin{enumerate}
\item For $f\in\omega^\omega$ and $\varphi\in\cS$,
    if $f(n)\in\varphi(n)$ holds
    for all but finitely many $n<\omega$,
    then $H_f\subseteq H_\varphi$.
\item For $\varphi,\psi\in\cS$,
    if $\psi(n)\subseteq\varphi(n)$ holds
    for all but finitely many $n<\omega$,
    then $H_\psi\subseteq H_\varphi$.
\end{enumerate}
Note that the reversed implications in the above statements do not
hold in general.

We will now describe a canonical procedure for constructing a
nonempty closed set disjoint from $H_\varphi$. For a slalom
$\varphi\in\cS$, define a function $r_\varphi\in\omega^\omega$ by
induction on $n<\omega$ as follows: $r_\varphi(0)=0$, and for
$1\leq n<\omega$, let
\[
r_\varphi(n)=\min\{i<\omega\st
    C^n_i\subseteq C^{n-1}_{r_\varphi(n-1)}
        \ssm\bigcup_{j\in\varphi(n)}C^n_j\}.
\]
This induction goes well because, by the choice of $h$, we have
$\mu(C^{n-1}_k)\geq(n+1)\cdot\mu(C^n_j)$ for $j,k<\omega$.

Let $R_\varphi=\bigcap_{n<\omega}C^n_{r_\varphi(n)}$. $R_\varphi$
is a nonempty closed set, because it is the intersection of a
decreasing sequence of nonempty closed sets in a compact space.
Let $A_\varphi=\bigcup_{n<\omega}\bigcup_{i\in\varphi(n)}C^n_i$.
Then clearly $H_\varphi\subseteq A_\varphi$. By the construction
of $r_\varphi$, we have $R_\varphi\cap A_\varphi=\emptyset$, and
hence $R_\varphi\cap H_\varphi=\emptyset$.

For $\varphi,\psi\in\cS$, if $r_\varphi(n)\in\psi(n)$ for
infinitely many $n<\omega$, then $R_\varphi\subseteq H_\psi$ and
hence $H_\psi\not\subseteq H_\varphi$.

\section{Localization forcing}\label{sec:loc}

In this section, we will introduce a modified form of
\emph{localization forcing} $\LOC$, which is defined in
\cite[Section~3.1]{BaJ:set}.

Let $\cT=\bigcup_{n<\omega}\prod_{i<n}[\omega]^{\leq i}$. A
condition $p$ of $\LOC$ is of the form $p=(s^p,F^p)$, where
$s^p\in\cT$, $F^p\subseteq\omega^\omega$ and
$\size{F^p}\leq\lh{s^p}$. For conditions $p,q$ in $\LOC$, $p\leq
q$ if $s^p\supseteq s^q$, $F^p\supseteq F^q$, and for each
$n\in\lh{s^p}\ssm\lh{s^q}$ and $f\in F^q$ we have $f(n)\in
s^p(n)$.

It is easy to see the following.
\begin{enumerate}
\item For each $n<\omega$,
    the set $\{q\in\LOC\st\lh{s^q}\geq n\}$ is dense in $\LOC$.
\item For each $f\in\omega^\omega$,
    the set $\{q\in\LOC\st f\in F^q\}$ is dense in $\LOC$.
\item $\LOC$ is $\sigma$-linked,
and hence it satisfies ccc.
\end{enumerate}

Let $\V$ be a ground model, and $G$ a $\LOC$-generic filter over
$\V$. In $\V[G]$, let $\varphi_G=\bigcup\{s^p\st p\in G\}$. Then
$\varphi_G\in\cS$ and, for every $f\in\omega^\omega\cap\V$, for
all but finitely many $n<\omega$ we have $f(n)\in\varphi_G(n)$.

Let $H_G=H_{\varphi_G}$. Then in $\V[G]$, by the observation in
Section~\ref{sec:combnull}, for every Borel null set $X\subseteq
2^\omega$ which is coded in $\V$, we have $X\subseteq H_G$.

Now we define a modified form of localization forcing.

\begin{definition}
Define $\LOC^*$ as follows. A condition $p$ of $\LOC^*$ is of the
form $p=(s^p,w^p,F^p)$, where
\begin{enumerate}
\item $s^p\in\cT$,
    $w^p<\omega$,
    $F^p\subseteq\omega^\omega$, and
\item $\size{F^p}\leq w^p\leq\lh{s^p}$.
\end{enumerate}
For $p,q\in\LOC^*$, $p\leq q$ if
\begin{enumerate}
\setcounter{enumi}{2}
\item
    $s^q\subseteq s^p$,
    $w^q\leq w^p$,
    $F^q\subseteq F^p$,
    and
    for $n\in\lh{s^p}\ssm\lh{s^q}$ and $f\in F^q$
    we have $f(n)\in s^p(n)$;
\item $w^p\leq w^q+(\lh{s^p}-\lh{s^q})$;
\item For $n\in\lh{s^p}\ssm\lh{s^q}$,
    we have $\size{s^p(n)}\leq w^q+(n-\lh{s^q})$.
\end{enumerate}
\end{definition}

We show that the forcing $\LOC^*$ has similar properties to
$\LOC$.

\begin{lemma}\label{lem:locprolong}
For each $n<\omega$, the set $\{q\in\LOC^*\st\lh{s^q}\geq n\}$ is
dense in $\LOC^*$.
\end{lemma}

\begin{proof}
Easy.
\end{proof}

\begin{lemma}\label{lem:locaddfunc}
For each $f\in\omega^\omega$, the set $\{q\in\LOC^*\st f\in F^q\}$
is dense in $\LOC^*$.
\end{lemma}

\begin{proof}
Fix $p\in\LOC^*$ and $f\in\omega^\omega$. Define $q=(s^q,w^q,F^q)$
as follows: $\lh{s^q}=\lh{s^p}+1$, $s^q\restrictedto\lh{s^p}=s^p$,
$s^q(\lh{s^p})=\{f(\lh{s^p})\st f\in F^p\}$, $w^q=w^p+1$ and
$F^q=F^p\cup\{f\}$. It is easy to see that $q\in\LOC^*$ and $q\leq
p$.
\end{proof}

\begin{lemma}
$\LOC^*$ is $\sigma$-linked, and hence it satisfies ccc.
\end{lemma}

\begin{proof}
It is easily seen that the set $L=\{p\in\LOC^*\st w^p\geq
2\cdot\size{F^p}\}$ is dense in $\LOC^*$. For each $s\in\cT$ and
$w\leq\lh{s}$, let $L_{s,w}=\{p\in L\st s^p=s\text{ and }w^p=w\}$.
Then $L=\tbigcup\{L_{s,w}\st s\in\cT\text{ and }w\leq\lh{s}\}$
and, for each $s\in\cT$ and $w\leq\lh{s}$, any two conditions in
$L_{s,w}$ are compatible.
\end{proof}

Let $\V$ be a ground model, and $G$ a $\LOC^*$-generic filter over
$\V$. In $\V[G]$, let $\varphi_G=\bigcup\{s^p\st p\in G\}$. Then,
by Lemmata~\ref{lem:locprolong} and \ref{lem:locaddfunc}, we have
$\varphi_G\in\cS$ and, for every $f\in\omega^\omega\cap\V$, for
all but finitely many $n<\omega$ we have $f(n)\in\varphi_G(n)$.

Let $H_G=H_{\varphi_G}$. The following proposition follows from
the observation in Section~\ref{sec:combnull}.

\begin{proposition}\label{prop:locamoeba}
Let $\V$ be a ground model and $G$ a $\LOC^*$-generic filter over
$\V$. Then in $\V[G]$, for every Borel null set $X\subseteq
2^\omega$ which is coded in $\V$, we have $X\subseteq H_G$.
\end{proposition}

As we observed in Section~\ref{sec:combnull}, in $\V[G]$, we can
define $r_{\varphi_G}$ and $R_{\varphi_G}$ from $\varphi_G$. Note
that, in this context, every $x\in R_{\varphi_G}$ is a random real
over $\V$. We can naturally define a $\LOC^*$-name $\dot{r}$ for
$r_{\varphi_G}$ so that, for $p\in\LOC^*$, if $\lh{s^p}=n$ then
$p$ decides the value of $\dot{r}\restrictedto n$, because
$r_{\varphi_G}\restrictedto n$ depends only on
$\varphi_G\restrictedto n$.

\section{Hechler's theorem for the null ideal}\label{sec:wfi}

In this section, we will construct a ccc forcing notion which
yields Hechler's theorem for the null ideal. The idea is to use
localization forcing at each step, instead of the dominating real
partial order used in Hechler's construction.

Let $(Q,\leq)$ be a partially ordered set such that every
countable subset of $Q$ has a strict upper bound in $Q$, that is,
for every countable set $A\subseteq Q$ there is $b\in Q$ such that
$a<b$ for all $a\in A$. Extend the order to $Q^*=Q\cup\{Q\}$ by
letting $a<Q$ for all $a\in Q$.

Fix a well-founded cofinal subset $R$ of $Q$. Define the rank
function on the well-founded set $R^*=R\cup\{Q\}$ in the usual
way. For $a\in Q\ssm R$, let $\rk{a}=\min\{\rk{b}\st b\in
R^*\text{ and }a<b\}$. For $x,y\in Q^*$, we say $x\ll y$ if $x<y$
and $\rk{x}<\rk{y}$. For $x\in Q^*$, let $Q_x=\{y\in Q\st y\ll
x\}$.

For $D\subseteq Q$ and $\xi\leq\rk{Q}$, let $D_{<\xi}=\{y\in
D\st\rk{y}<\xi\}$, $D_\xi=\{y\in D\st\rk{y}=\xi\}$, and for $x\in
Q$ with $\rk{x}=\xi$, let $D_{\leq x}=\{y\in D_\xi\st y\leq x\}$.

For $D\subseteq Q$, let $\bar{D}=\{\rk{x}\st x\in D\}$.

For $E\subseteq D\subseteq Q$, we say \emph{$E$ is downward closed
in }$D$ if, for $x\in E$ and $y\in D$ if $y\leq x$ then $y\in E$.
When $E$ is downward closed in $Q$, we simply say \emph{$E$ is
downward closed}.

\begin{definition}\label{def:iteration}
We define forcing notions $\bN_a$ for $a\in Q^*$ by induction on
$\rk{a}$. For $a\in Q^*$, the conditions $p$ of $\bN_a$ are all
objects of the form $p=\{(s^p_x,w^p_x,F^p_x)\st x\in D^p\}$ which
satisfy the following properties.
\begin{enumerate}
\item\label{item:cdomain}
    $D^p$ is a finite subset of $Q_a$;
\item\label{item:catom}
    For $x\in D^p$,
    $s^p_x\in\cT$,
    $w^p_x<\omega$,
    $F^p_x$ is a finite set of $\bN_x$-names
    for functions in $\omega^\omega$,
    and
    $\size{F^p_x}\leq w^p_x$;
\item\label{item:cwtotal}
    For $x\in D^p$,
    $\tsum\{w^p_z\st z\in D^p_{\leq x}\}\leq\lh{s^p_x}$;
\item\label{item:clength}
    For $x,y\in D^p$,
if $\rk{x}=\rk{y}$
    then $\lh{s^p_x}=\lh{s^p_y}$.
\end{enumerate}
As in the definition of iterated forcing, it is necessary to limit
the collection of names in clause~2 so that $\bN_a$ is not a
proper class. We leave it understood that by a name for an element
of $\omega^\omega$ is meant a nice name for a subset of
$(\omega\times\omega)\check{}$ in the sense \cite[VII
5.11]{Ku:set} which is forced by the weakest condition to name an
element of $\omega^\omega$.

Throughout this paper, for a condition $p$ in $\bN_a$, we always
use the notation $D^p$, $s^p_x$, $w^p_x$ and $F^p_x$ to denote
respective components of $p$. Also, for $p\in\bN_a$ and
$\xi\in\bar{D^p}$, let $l^p_\xi$ be the length of $s^p_x$ for
$x\in D^p_\xi$.

For $p\in\bN_a$ and $b\in Q_a$, define $p\restrictedto b\in\bN_b$
by letting $p\restrictedto b=\{(s^p_x,w^p_x,F^p_x)\st x\in D^p\cap
Q_b\}$.

For conditions $p,q$ in $\bN_a$, $p\leq q$ if:
\begin{enumerate}
\setcounter{enumi}{4}
\item\label{item:rdomain} $D^q\subseteq D^p$;
\item\label{item:ratom} For $x\in D^q$,
    $s^p_x\supseteq s^q_x$,
    $w^p_x\geq w^q_x$,
    $F^p_x\supseteq F^q_x$
    and,
    for all $n\in\lh{s^p_x}\ssm\lh{s^q_x}$ and $\dot{f}\in F^q_x$
    we have
    $p\restrictedto x\forces[\bN_x]{\dot{f}(n)\in s^p_x(n)}$;
\item\label{item:rcannibal}
    For $\xi\in\bar{D^q}$ and $x,y\in D^q_\xi$,
    if $x<y$,
    then
    for all $n\in l^p_\xi\ssm l^q_\xi$
    we have
    $s^p_x(n)\subseteq s^p_y(n)$;
\item\label{item:rwgrowth}
    For $\xi\in\bar{D^q}$,
    $\tsum\{w^p_x\st x\in D^p_\xi\}
        \leq\tsum\{w^q_x\st x\in D^q_\xi\}+(l^p_\xi-l^q_\xi)$;
\item\label{item:rsgrowth}
    For $\xi\in\bar{D^q}$,
    $E\subseteq D^q_\xi$ which is downward closed in $D^q_\xi$
    and $n\in l^p_\xi\ssm l^q_\xi$,
    we have
    \[\size{\tbigcup\{s^p_x(n)\st x\in E\}}
        \leq\tsum\{w^q_x\st x\in E\}+(n-l^q_\xi).\]
\end{enumerate}
\end{definition}

\begin{remark}\label{rem:rwgrowthdiscardterms}
If $p\leq q$, then for any $\xi\in\bar{D^q}$ and $E\subseteq
D^p_\xi$ we can discard the terms with indices not in $E$ from
both sides of the inequality in clause~\ref{item:rwgrowth} (using
$w^p_x\geq w^q_x$ from clause~\ref{item:ratom}) to get
\[
\tsum\{w^p_x\st x\in E\}
    \leq\tsum\{w^q_x\st x\in E\cap D^q_\xi\}+(l^p_\xi-l^q_\xi).
\]
\end{remark}

We now verify that Definition~\ref{def:iteration} does indeed
define a partial order. (Reflexivity is clear, but we need to
prove transitivity.) The simple observation in part~(c) of the
following proposition justifies not mentioning $a$ in the notation
$\leq$ for the order relation on $\bN_a$.

\begin{proposition}\label{transitivity}
We have the following properties.
\begin{enumerate}
\item[{\rm(a)}]
    For any conditions $p,q\in \bN_a$,
    if $p\leq q$ then for any $b\in Q_a$,
    $p\restrictedto b\leq q\restrictedto b$.
\item[{\rm(b)}]
    The order relation on $\bN_a$ is transitive.
\item[{\rm(c)}]
    For any $a,b\in Q^*$,
    if $p,q\in\bN_a\cap \bN_b$,
    then $p\leq q$
    in $\bN_a$ if and only if $p\leq q$ in $\bN_b$.
\end{enumerate}
\end{proposition}

\begin{proof}
(a) and (b) are proven simultaneously by induction on the rank of
$a$. Note that part (b) of the induction hypothesis ensures that
for $p,q\in\bN_a$ and $x\in D^q\subseteq Q_a$, $\bN_x$ is a
well-defined partial order and hence the last part of clause~6
makes sense.

(a) All but the last part of clause~\ref{item:ratom} and
clause~\ref{item:rwgrowth} in the definition of $p\restrictedto
b\leq q\restrictedto b$ are inherited directly from the
corresponding clauses for $p\leq q$. The last part of
clause~\ref{item:ratom} holds because for $x\in D^{q\restrictedto
b}=D^q\cap Q_b$, $(p\restrictedto b)\restrictedto x=p\restrictedto
x$. There remains to check clause~\ref{item:rwgrowth}. Let $\xi\in
\bar D^{q\restrictedto b}$. Using clause~\ref{item:rwgrowth} for
$p\leq q$ and the fact that $w^p_x\geq w^q_x$ whenever both are
defined, we have
\begin{align*}
\tsum\{w^{p\restrictedto b}_x\st x\in D^{p\restrictedto b}_\xi\}
 & = \tsum\{w^p_x\st x\in D^{p\restrictedto b}_\xi\}\\
 & = \tsum\{w^p_x\st x\in D^p_\xi\}
    - \tsum\{w^p_x\st x\in D^p_\xi\ssm Q_b\}\\
 & \leq  \tsum\{w^q_x\st x\in D^q_\xi\} + (l^p_\xi - l^q_\xi)
    - \tsum\{w^p_x\st x\in D^p_\xi\ssm Q_b\}\\
 & \leq  \tsum\{w^q_x\st x\in D^q_\xi\} + (l^p_\xi - l^q_\xi)
    - \tsum\{w^q_x\st x\in D^q_\xi\ssm Q_b\}\\
 & =  \tsum\{w^{q\restrictedto b}_x\st x\in D^{q\restrictedto b}_\xi\}
    + (l^{p\restrictedto b}_\xi - l^{q\restrictedto b}_\xi).
\end{align*}

(b) Suppose that $a\in Q^*$, $p,q,r\in\bN_a$ and $p\leq q\leq r$.
We must show $p\leq r$.

For the last part of clause~\ref{item:ratom}, suppose we have
$x\in D^r_\gamma$, $n\in l^p_\gamma\ssm l^r_\gamma$, $\dot f\in
F^r_x$. If $n\in l^p_\gamma\ssm l^q_\gamma$, then because $\dot
f\in F^r_x\subseteq F^q_x$, the fact that $p\leq q$ gives
$p\restrictedto x\Vdash_{\bN_x}\dot f(n)\in s^p_x(n)$. If $n\in
l^q_\gamma\ssm l^r_\gamma$, then the fact that $q\leq r$ gives
$q\restrictedto x\Vdash_{\bN_x}\dot f(n)\in s^q_x(n)$. We have
$s^p_x(n)=s^q_x(n)$ by the first part of clause~\ref{item:ratom}
for $p\leq q$. Also, $p\restrictedto x\leq q\restrictedto x$ by
part~(a). Thus, $p\restrictedto x\Vdash_{\bN_x}\dot f(n)\in
s^p_x(n)$. We now check clause~\ref{item:rsgrowth} and leave the
other clauses for the reader. Fix $\xi\in\bar{D^r}$, $E\subseteq
D^r_\xi$ which is downward closed in $D^r_\xi$ and $n\in
l^p_\xi\ssm l^r_\xi$. Let $E^q$ be the downward closure of $E$ in
$D^q_\xi$. If $n\in l^q_\xi\ssm l^r_\xi$, then
\begin{align*}
    \size{\tbigcup\{s^p_x(n)\st x\in E\}}
    & = \size{\tbigcup\{s^q_x(n)\st x\in E\}}\\
    & \leq \tsum\{w^r_x\st x\in E\}+(n-l^r_\xi)
\end{align*}
because of clause~\ref{item:rsgrowth} for $q\leq r$. If $n\in
l^p_\xi\ssm l^q_\xi$, then
\begin{align*}
    \size{\tbigcup\{s^p_x(n)\st x\in E\}}
    & \leq  \size{\tbigcup\{s^p_x(n)\st x\in E^q\}}   \\
    & \leq  \tsum\{w^q_x\st x\in E^q\}+(n-l^q_\xi)    \\
    & \leq  \tsum\{w^r_x\st x\in E\}
            +(l^q_\xi-l^r_\xi)+(n-l^q_\xi)  \\
    & = \tsum\{w^r_x\st x\in E\}
            +(n-l^r_\xi).
\end{align*}
The second inequality follows from clause~\ref{item:rsgrowth} for
$p\leq q$ and the third from Remark~\ref{rem:rwgrowthdiscardterms}
for $q\leq r$. Hence we have $p\leq r$.

(c) The definition of the order on $\bN_a$ makes no mention of
$a$.
\end{proof}

\begin{definition}
For a downward closed set $A\subseteq Q$, let
$\bN_A=\{p\in\bN_Q\st D^p\subseteq A\}$, and for $p\in\bN_Q$, we
define $p\restrictedto A\in\bN_A$ by letting $p\restrictedto
A=\{(s^p_x,w^p_x,F^p_x)\st x\in D^p\cap A\}$. For $\xi\leq\rk{Q}$,
let $\bN_\xi=\bN_{Q_{<\xi}}$ and $p\restrictedto\xi=p\restrictedto
Q_{<\xi}$. Also, for $\xi\leq\rk{Q}$, let $p\restrictedto\{\xi\}
    =\{(s^p_x,w^p_x,F^p_x)\st x\in D^p_\xi\}\in\bN_{\xi+1}$
and $p\restrictedto [\xi,\infty)
    =\{(s^p_x,w^p_x,F^p_x)\st x\in D^p\ssm Q_{<\xi}\}\in\bN_Q$.
\end{definition}

In this notation, $\bN_a=\bN_{Q_a}$ for $a\in Q$, and $\bN_Q$ has
the same meaning if we consider the subscript $Q$ either as an
element of $Q^*$ or as a subset of $Q$.

Clearly $A\subseteq B\subseteq Q$ implies $\bN_A\subseteq
\bN_B\subseteq \bN_Q$. We are going to prove that, if $A\subseteq
B$, then $\bN_A$ is completely embedded into $\bN_B$. This is a
fundamental principle of iterated forcing.

The following lemma, which is a special case of this principle, is
easily checked.

\begin{lemma}\label{lem:easycompleteembedding}
If $B$ is a downward closed subset of $Q$, $\xi\leq\rk{Q}$,
$p\in\bN_B$ and $q\in \bN_{B_{<\xi}}$ extends $p\restrictedto\xi$,
then $q\cup p\restrictedto[\xi,\infty)$ belongs to $\bN_B$ and
extends both $p$ and $q$. In particular, $\bN_{B_{<\xi}}$ is
completely embedded into $\bN_B$.
\end{lemma}

Using this lemma, we prove the following.

\begin{lemma}\label{lem:completeembedding}
For downward closed sets $A,B\subseteq Q$, if $A\subseteq B$, then
$\bN_A$ is completely embedded into $\bN_B$ by the identity map.
\end{lemma}

\begin{proof}
It is easy to see that the compatibility of conditions in $\bN_A$
is the same either in $\bN_A$ or in $\bN_B$. We show that, for
$p\in\bN_B$ and $r\in\bN_A$, if $r\leq p\restrictedto A$ then
there is $q\in\bN_B$ satisfying $q\leq p$ and $q\leq r$. We will
proceed by induction on $\sup{\bar{A}}$.

Suppose that $p\in\bN_B$, $r\in\bN_A$ and $r\leq p\restrictedto
A$. Let $\gamma=\max\bar{D^r}$. By the induction hypothesis, there
is $q_{<\gamma}\in\bN_{B_{<\gamma}}$ satisfying $q_{<\gamma}\leq
p\restrictedto\gamma$ and $q_{<\gamma}\leq r\restrictedto\gamma$.

For $x\in D^r_\gamma$, let $(s_x,w_x,F_x)=(s^r_x,w^r_x,F^r_x)$.
For $x\in D^p_\gamma\ssm D^r_\gamma$, let
$(s_x,w_x,F_x)=(s^p_x,w^p_x,F^p_x)$.

Let
\[
L=\tsum\{w_x\st  x\in D^p_\gamma\cup D^r_\gamma\}+l^r_\gamma.
\]
By the induction hypothesis, for each $x\in D^p_\gamma\cup
D^r_\gamma$, $\bN_x$ is completely embedded into
$\bN_{B_{<\gamma}}$ and so each $\dot{f}\in F_x$ is an
$\bN_{B_{<\gamma}}$-name. Choose $q^*\in\bN_{B_{<\gamma}}$ so that
$q^*\leq q_{<\gamma}$ and $q^*$ decides the values of
$\dot{f}\restrictedto L$ for all $\dot{f}\in\tbigcup\{F_x\st  x\in
D^p_\gamma\cup D^r_\gamma\}$. For $x\in D^p_\gamma\cup D^r_\gamma$
and $n\in L\ssm\lh{s_x}$, let $K_{x,n}\subseteq\omega$ be the set
satisfying $q^*\forces{K_{x,n}=\{\dot{f}(n)\st \dot{f}\in F_x\}}$.

Define $s^*_x$ for $x\in D^p_\gamma\cup D^r_\gamma$ in the
following way: If $x\in D^r_\gamma$, then $\lh{s^*_x}=L$,
$s^*_x\restrictedto l^r_\gamma=s_x$, and for $n\in L\ssm
l^r_\gamma$,
\[
    s^*_x(n)=\tbigcup\{K_{z,n}\st  z\in D^r_{\leq x}\}.
\]
If $x\in D^p_\gamma\ssm D^r_\gamma$, then $\lh{s^*_x}=L$,
$s^*_x\restrictedto l^p_\gamma=s_x$, and for $n\in L\ssm
l^p_\gamma$,
\[
    s^*_x(n)=
    \begin{cases}

    \tbigcup\{s_z(n)\st  z\in D^p_{\leq x}\cap D^r_\gamma\}
       \cup\tbigcup\{K_{z,n}\st  z\in D^p_{\leq x}\ssm D^r_\gamma\}
       &   \text{if }l^p_\gamma\leq n<l^r_\gamma  \\
    \tbigcup\{K_{z,n}\st  z\in (D^p_\gamma\cup D^r_\gamma)_{\leq x}\}
       &   \text{if }l^r_\gamma\leq n<L,\
           \gamma\in \bar D^{p\restrictedto A}\\
    \tbigcup\{K_{z,n}\st z\in D^p_{\leq x}\}
       &   \text{if }l^r_\gamma\leq n<L,\
           \gamma\not\in \bar D^{p\restrictedto A}
    \end{cases}
\]
Now we define $q=\{(s^q_x,w^q_x,F^q_x)\st  x\in D^q\}$ by the
following:
\begin{enumerate}
\item $D^q=D^p\cup D^{q^*}\cup D^r_\gamma$;
\item
    For $x\in D^{q^*}$,
    $(s^q_x,w^q_x,F^q_x) = (s^{q^*}_x,w^{q^*}_x,F^{q^*}_x)$;
\item
    For $x\in D^p_\gamma\cup D^r_\gamma$,
    $(s^q_x,w^q_x,F^q_x) = (s^*_x,w_x,F_x)$;
\item
    For $x\in D^p\ssm Q_{<\gamma+1}$,
    $(s^q_x,w^q_x,F^q_x) = (s^p_x,w^p_x,F^p_x)$.
\end{enumerate}
We now check that $q\in\bN_B$. The conditions of
Definition~\ref{def:iteration} are satisfied below (resp.\ above)
rank $\gamma$ because $q^*$ (resp.\ $p$) is a condition. Consider
what they say at rank $\gamma$. The first clause is trivial. The
fourth holds because the $s^{q}_x$'s all have domain $L$. The
third clause can be checked in two cases.
\begin{enumerate}
\item[(i)]
If $x\in D^{r}_\gamma$, then $D^{q}_{\leq x}=(D^p\cup D^{r})_{\leq
x}=D^{r}_{\leq x}$, so $\tsum\{w^{q}_z\st z\in D^{q}_{\leq x}\}
    =\tsum\{w^{r}_z\st z\in D^{r}_{\leq x}\}
    \leq l^{r}_\gamma\leq L$.
\item[(ii)]
If $x\in D^{p}_\gamma\ssm D^{r}_\gamma$, then $D^{q}_{\leq
x}=D^{p}_{\leq x}\cup D^{r}_{\leq x}$, so $ \tsum\{w^{q}_z\st z\in
D^{q}_{\leq x}\}
    =\tsum\{w_z\st z\in D^{p}_{\leq x}\cup D^{r}_{\leq x}\}
    \leq L.$
\end{enumerate}
For the second, all the requirements except that the $s^{q}_x$'s
are partial slaloms follow from the fact that $p$ and $r$ are
conditions. We need to check that $\size{s^*_x(n)}\leq n$ for each
relevant $n$. If $x\in D^{r}_\gamma$, then for $l^{r}_\gamma\leq
n<L$, we have $\size{s^*_x(n)}
    \leq\tsum\{w^{r}_z\st z\in D^{r}_{\leq x}\}
    \leq\lh{s^{r}_x}=l^{r}_\gamma\leq n$.
If $x\in D^p_\gamma\ssm D^{r}_\gamma$, we consider four cases.

Case 1. $l^p_\gamma\leq n<l^{r}_\gamma$ and $\gamma\in \bar
D^{p\restrictedto A}$. Definition~\ref{def:iteration}(9) for
$r\leq p\restrictedto A$ with $E=D^p_{\leq x}\cap D^{r}_\gamma$
gives
\begin{align*}
\size{s^*_x(n)} & \leq  \tsum\{w^{p}_z\st z\in E\} +
(n-l^p_\gamma) +
    \tsum\{w^{p}_z\st z\in D^{p}_{\leq x}\ssm E\}\\
           & =  \tsum\{w^{p}_z\st z\in D^{p}_{\leq x}\} +
           (n-l^p_\gamma)\\
           & \leq  l^p_\gamma + (n - l^p_\gamma) = n.
\end{align*}

Case 2. $l^p_\gamma\leq n<l^{r}_\gamma$ and $\gamma\not\in \bar
D^{p\restrictedto A}$. In this case, $D^p_{\leq x}\cap
D^r_\gamma\subseteq D^p_\gamma\cap A=\emptyset$, so
$\size{s^*_x(n)}
    \leq \tsum\{w^{p}_z\st z\in D^{p}_{\leq x}\}
    \leq l^p_\gamma\leq n$.

Case 3. $l^{r}_\gamma\leq n<L$ and $\gamma\in \bar
D^{p\restrictedto A}$. Definition~\ref{def:iteration}(8) for
$r\leq p\restrictedto A$ gives $\tsum\{w^{r}_z\st z\in
D^{r}_\gamma\}
    \leq \tsum\{w^p_z\st z\in D^{p\restrictedto A}_\gamma\}
    + (l^{r}_\gamma-l^p_\gamma)$.
Removing terms with $z\not\leq x$ from both sides (see
Remark~\ref{rem:rwgrowthdiscardterms}) gives
\[
\tsum\{w^{r}_z\st z\in D^{r}_{\leq x}\}
    \leq \tsum\{w^p_z\st z\in D^p_{\leq x}\cap A\}
    +(l^{r}_\gamma-l^p_\gamma).
\]
From the formula for $s^*_x(n)$ we now get
\begin{align*}
\size{s^*_x(n)}
    & \leq  \tsum\{w^{r}_z\st z\in D^{r}_{\leq x}\}+
        \tsum\{w^{p}_z\st z\in D^p_{\leq x}\ssm A\}\\
    & \leq  \tsum\{w^p_z\st z\in D^p_{\leq x}\cap A\} +
        (l^{r}_\gamma-l^p_\gamma)
        + \tsum\{w^{p}_z\st z\in D^p_{\leq x}\ssm A\}\\
    & =  \tsum\{w^p_z\st z\in D^p_{\leq x}\} +
        (l^{r}_\gamma-l^p_\gamma)\\
    & \leq  l^p_\gamma + (l^{r}_\gamma - l^p_\gamma)
        = l^{r}_\gamma\leq n.
\end{align*}

Case 4. $l^{r}_\gamma\leq n<L$ and $\gamma\not\in \bar
D^{p\restrictedto A}$. In this case we have $\size{s^*_x(n)}
    \leq \tsum\{w^{p}_z\st z\in D^{p}_{\leq x}\}\leq l^p_\gamma\leq n$.

Thus, $q$ is a condition.

We now check Definition~\ref{def:iteration}(5--9) for $q\leq r$
and $q\leq p$. Clause~5 follows from the definition of $q$. For
clauses~6--9, first note that below rank $\gamma$, they hold
because $q^*\leq p\restrictedto\gamma$ and $q^*\leq
r\restrictedto\gamma$. Consider what happens at rank $\gamma$.
Clause~6 holds because for $x\in D^{p}_\gamma\cup D^{r}_\gamma$
and all the relevant values of $\dot f$ and $n$, we have from the
definitions that $q^*\Vdash\dot f(n)\in K_{x,n}$ and
$K_{x,n}\subseteq s^*_x(n)$. For clause~7, we consider three
cases. Let $x<y$ be elements of $D^p_\gamma\cup D^{r}_\gamma$.
\begin{enumerate}
\item[(i)]
If $x,y\in D^{r}_\gamma$, then for checking $q\leq r$, just use
the monotonicity of $s^*_x(n)$ as a function of $x$. For checking
$q\leq p$ (so now we assume $x,y\in D^p_\gamma$ as well), we also
need to consider values of $n$ such that $l^p_\gamma\leq
n<l^{r}_\gamma$. But then $s^*_x(n)=s^{r}_x(n)\subseteq s^{r}_y(n)
= s^*_y(n)$ because $r\leq p\restrictedto A$.

This is the only case to consider for checking clause~7 for $q\leq
r$ at stage $\gamma$. The remaining cases deal with checking
$q\leq p$. Note that if $y\in D^{r}_\gamma\cap
D^p_\gamma=D^{p}_\gamma\cap A$ then also $x \in D^{r}_\gamma\cap
D^p_\gamma$ since $A$ is downward closed.
\item[(ii)]
If $x,y\in D^p_\gamma\ssm D^{r}_\gamma$, use the monotonicity of
$s^*_x(n)$ as a function of $x$.
\item[(iii)]
If $x\in D^{r}_\gamma\cap D^p_\gamma$ and $y\in D^p_\gamma\ssm
D^{r}_\gamma$, then consider first a value of $n$ such that
$l^p_\gamma\leq n<l^{r}_\gamma$. We have $s^*_x(n)=s_x(n)
    \subseteq \tbigcup\{s_z(n)\st z\in D^p_{\leq y}\cap D^{r}_\gamma\}
    \subseteq s^*_y(n)$.
Next consider $n$ such that $l^{r}_\gamma\leq n<L$. We have
$s^*_x(n)
    = \tbigcup\{K_{z,n}\st  z\in D^{r}_{\leq x}\}
    \subseteq \tbigcup\{K_{z,n}
        \st z\in (D^p_\gamma\cup D^{r}_\gamma)_{\leq y}\}
    =s^*_y(n)$.
\end{enumerate}
This takes care of clause~7. Clause~8 follows from the fact that
from the definition of $L$ we have $\tsum\{w_x\st x\in
D^{r}_\gamma\cup D^p_\gamma\}
    \leq L-l^{r}_\gamma\leq L-l^{p}_\gamma$.
For clause 9, first we check $q\leq r$. If $E\subseteq
D^{r}_\gamma$ is downward closed in $D^{r}_\gamma$ and
$l^{r}_\gamma\leq n<L$, then $\size{\tbigcup\{s^*_x(n)\st x\in
E\}}
    =\size{\tbigcup\{K_{x,n}\st x\in E\}}
    \leq\tsum\{w^{r}_x\st x\in E\}$.
Next we check $q\leq p$. Suppose $\gamma\in\bar D^p$ and let
$E\subseteq D^p_\gamma$ be downward closed. Consider four cases.

Case 1. $l^p_\gamma\leq n<l^{r}_\gamma$ and $\gamma\in \bar
D^{p\restrictedto A}$. Using Definition~\ref{def:iteration}(9) for
$r\leq p\restrictedto A$ and the fact that $E\cap A$ is downward
closed in $D^{p\restrictedto A}$, we have
\begin{align*}
\size{\tbigcup\{s^*_x(n)\st x\in E\}}
    &= \size{\tbigcup\{s^*_x(n)\st x\in E\cap A\}
        \cup\tbigcup\{s^*_x(n)\st x\in E\ssm A\}}   \\
    &= \size{\tbigcup\{s^{r}_x(n)\st x\in E\cap A\}
        \cup \tbigcup\{K_{x,n}\st x\in E\ssm A\}}   \\
    &\leq \tsum\{w^{p\restrictedto A}_x\st x\in E\cap A\}+(n-l^p_\gamma) +
        \tsum\{w^p_x\st x\in E\ssm A\}  \\
    &= \tsum\{w^{p}_x\st x\in E\}+(n-l^p_\gamma).
\end{align*}

Case 2. $l^p_\gamma\leq n<l^{r}_\gamma$ and $\gamma\not\in \bar
D^{p\restrictedto A}$. Then $E\cap A=\emptyset$, and the
calculation for case~1 reduces to
\begin{align*}
\size{\tbigcup\{s^*_x(n)\st x\in E\}}
    &= \size{\tbigcup\{s^*_x(n)\st x\in E\ssm A\}}  \\
    &= \size{\tbigcup\{K_{x,n}\st x\in E\ssm A\}}   \\
    &\leq \tsum\{w^p_x\st x\in E\ssm A\}    \\
    &\leq \tsum\{w^{p}_x\st x\in E\}+(n-l^p_\gamma).
\end{align*}

Case 3. $l^{r}_\gamma\leq n<L$ and $\gamma\in \bar
D^{p\restrictedto A}$. Let $E^r$ be the downward closure in
$D^r_\gamma$ of $E\cap A=E\cap D^{p\restrictedto A}_\gamma$. Using
Definition~\ref{def:iteration}(8) for $r\leq p\restrictedto A$ and
removing terms with $z\not\in E^r$ from both sides gives
\[
\tsum\{w^{r}_z\st z\in E^r\}\leq \tsum\{w^p_z\st z\in E\cap A\}
    +(l^{r}_\gamma-l^p_\gamma).
\]
Then we get
\begin{align*}
\size{\tbigcup\{s^*_x(n)\st x\in E\}}
    &= \size{\tbigcup\{K_{z,n}\st z\in E^r\}
        \cup\{K_{z,n}\st z\in E\ssm D^r_\gamma\}}   \\
    &\leq \tsum\{w^r_z\st z\in E^r\}+ \tsum\{w^p_z\st z\in E\ssm A\}    \\
    &\leq \tsum\{w^{p}_z\st z\in E\cap A\}
        + ( n - l^p_\gamma )
        + \tsum\{w^p_z\st z\in E\ssm A\}    \\
    &\leq \tsum\{w^{p}_z\st z\in E\} + ( n - l^p_\gamma ).
\end{align*}

Case 4. $l^{r}_\gamma\leq n<L$ and $\gamma\not\in \bar
D^{p\restrictedto A}$. We have
\[
\size{\tbigcup\{s^*_x(n)\st x\in E\}}
    = \size{\tbigcup\{K_{z,n}\st z\in E\}}
    \leq \tsum\{w^{p}_z\st z\in E\}.
\]

Thus, $q\leq r$. The proof that $q\leq p$ is completed by
appealing to Lemma~\ref{lem:easycompleteembedding}.
\end{proof}

The following definition and lemma provide a simple mechanism for
extending conditions.

\begin{definition}\label{def:precondition}
Let $B\subseteq Q$ be a downward closed set and
$\gamma\in\bar{B}$. $p'=\{(s^{p'}_x,w^{p'}_x,F^{p'}_x)\st x\in
D^{p'}\}$ is \emph{a $\gamma$-precondition of $\bN_B$} if $p'$
satisfies the following:
\begin{enumerate}
\item[\ref{item:cdomain}.]\label{item:pcdomain}
    $D^{p'}$ is a finite subset of $B$;
\item[\ref{item:catom}.]\label{item:pcatom}
For $x\in D^{p'}$, $s^{p'}_x\in\cT$, $w^{p'}_x<\omega$, $F^{p'}_x$
is a finite set of $\bN_x$-names for functions in $\omega^\omega$,
and $\size{F^{p'}_x}\leq w^{p'}_x$;
\item[\ref{item:cwtotal}$'$.]\label{item:pcwtotal}
    For $x\in D^{p'}\ssm D^{p'}_\gamma$,
    $\tsum\{w^{p'}_z\st z\in D^{p'}_{\leq x}\}\leq\lh{s^{p'}_x}$;
\item[\ref{item:clength}.]\label{item:pclength}
    For $x,y\in D^{p'}$,
if $\rk{x}=\rk{y}$
    then $\lh{s^{p'}_x}=\lh{s^{p'}_y}$.
\end{enumerate}
For $\xi\in\bar{D^{p'}}$, we will let $l^{p'}_\xi$ be the length
of $s^{p'}_x$ for $x\in D^{p'}_\xi$.

For $\gamma$-precondition $p'$ of $\bN_B$ and $p\in\bN_B$, we say
$p'$ is \emph{a $\gamma$-preextension of $p$} if
\begin{enumerate}
\item $D^{p'}\supseteq D^p$
    and $D^{p'}\ssm Q_{<\gamma+1}=D^p\ssm Q_{<\gamma+1}$;
\item $p'\restrictedto\gamma\leq p\restrictedto\gamma$;
\item For $x\in D^p_\gamma$,
    $s^{p'}_x=s^p_x$,
$F^{p'}_x=F^p_x$ and $w^{p'}_x\geq w^p_x$;
\item For $x\in D^{p'}_\gamma\ssm D^p_\gamma$,
    $F^{p'}_x=\emptyset$ and $w^{p'}_x=0$;
\item For $x\in D^p\ssm Q_{<\gamma+1}$,
    $(s^{p'}_x,w^{p'}_x,F^{p'}_x)=(s^p_x,w^p_x,F^p_x)$.
\end{enumerate}
\end{definition}

\begin{lemma}\label{lem:repair}
Let $B\subseteq Q$ be a downward closed set, $p\in\bN_B$,
$\gamma\in\bar{B}$, $p'=\{(s^{p'}_x,w^{p'}_x,F^{p'}_x)\st x\in
D^{p'}\}$ a $\gamma$-preextension of $p$ such that
$D^{p'}_\gamma\not=\emptyset$, and $N<\omega$. Then there is
$q\in\bN_B$ such that:
\begin{enumerate}
\item $q\leq p$ and $q\restrictedto\gamma\leq p'\restrictedto\gamma$;
\item $D^q_\gamma=D^{p'}_\gamma$ and,
    for $x\in D^q_\gamma$,
    $s^q_x\supseteq s^{p'}_x$,
    $w^q_x=w^{p'}_x$ and,
    $F^q_x=F^{p'}_x$;
\item $D^q\ssm Q_{<\gamma+1}=D^p\ssm Q_{<\gamma+1}$ and,
    for $x\in D^q\ssm Q_{<\gamma+1}$,
    $s^q_x=s^p_x$,
    $w^q_x=w^p_x$
    and $F^q_x=F^p_x$;
\item $l^q_\gamma\geq N$.
\end{enumerate}
\end{lemma}

\begin{proof}
Let $L=\max\{\tsum\{w^{p'}_x\st x\in D^{p'}_\gamma\} +
l^{p'}_\gamma,N\}$.

Note that clause~3 in the definition of ``$p'$ is a
$\gamma$-preextension of $p$'' ensures that
$l^{p'}_\gamma=l^p_\gamma$ as long as the latter is defined, i.e.,
as long as $\gamma\in \bar D^p$.

Using Lemma~\ref{lem:completeembedding}, choose
$q^*\in\bN_{B_{<\gamma}}$ so that $q^*\leq p'\restrictedto\gamma$
and $q^*$ decides the values of $\dot{f}\restrictedto L$ for all
$\dot{f}\in\tbigcup\{F^{p'}_x\st x\in D^{p'}_\gamma\}
    =\tbigcup\{F^p_x\st x\in D^p_\gamma\}$.
For $x\in D^{p}_\gamma$ and $n\in L\ssm l^{p'}_\gamma = L\ssm
l^{p}_\gamma$, let $K_{x,n}\subseteq\omega$ be the set satisfying
$q^*\forces{K_{x,n}=\{\dot{f}(n)\st\dot{f}\in F^{p}_x\}}$. Note
that $\size{K_{x,n}}\leq\size{F^p_x}\leq w^p_x$.

Define $s_x$ for $x\in D^{p'}_\gamma$ as follows: $\lh{s_x}=L$,
$s_x\restrictedto l^{p'}_\gamma=s^{p'}_x$, and for $n\in L\ssm
l^{p'}_\gamma$, if $x\in D^p_\gamma$ then
$s_x(n)=\tbigcup\{K_{z,n}\st z\in D^{p}_{\leq x}\}$ and if
$x\not\in D^p_\gamma$ then $s_x(n)=\emptyset$. Now we define
$q=\{(s^q_x,w^q_x,F^q_x)\st x\in D^q\}$ as follows:
\begin{enumerate}
\item $D^q=D^{q^*}\cup D^{p'}$;
\item
For $x\in D^{q^*}$,
$(s^q_x,w^q_x,F^q_x)=(s^{q^*}_x,w^{q^*}_x,F^{q^*}_x)$;
\item
For $x\in D^{p'}_\gamma$,
$(s^q_x,w^q_x,F^q_x)=(s_x,w^{p'}_x,F^{p'}_x)$;
\item
For $x\in D^q\ssm Q_{<\gamma+1}$,
$(s^q_x,w^q_x,F^q_x)=(s^{p'}_x,w^{p'}_x,F^{p'}_x)$.
\end{enumerate}
We now need to check that $q\in\bN_B$ and $q$ satisfies the
requirement. For $x\in D^{p'}_\gamma$, $l^{p'}_\gamma\leq n<L$, we
check that $\size{s_x(n)}\leq n$ and leave the rest of the
verification to the reader. If $x\not\in D^p_\gamma$, then
$s_x(n)=\emptyset$. Suppose now that $x\in D^p_\gamma$. Then
$\size{s_x(n)}= \size{\tbigcup\{K_{z,n}\st z\in D^p_{\leq x}\}}
    \leq \tsum\{w^p_z\st z\in D^p_{\leq x}\} \leq  \lh{s^p_x}
    = l^p_\gamma=l^{p'}_\gamma\leq n$.
\end{proof}

Next we prove that $\bN_Q$ satisfies ccc.

\begin{lemma}\label{lem:preccc}
Let $W$ be the collection of conditions $q\in\bN_Q$ satisfying the
following properties:
\begin{enumerate}
\item For all $x\in D^q$,
$2\cdot\size{F^q_x}\leq w^q_x$;
\item For all $\xi\in\bar{D^q}$,
    $2\cdot\tsum\{w^q_x\st x\in D^q_\xi\}\leq l^q_\xi$.
\end{enumerate}
Then $W$ is dense in $\bN_Q$.
\end{lemma}

\begin{proof}
By induction on $\xi\leq\rk{Q}$, we will show that $W_{<\xi}$ is
dense in $\bN_\xi$.

Fix $p\in\bN_\xi$ and let $\gamma=\max\bar{D^p}$. Define a
$\gamma$-preextension $p'$ of $p$ by the following: $D^{p'}=D^p$,
$p'\restrictedto\gamma=p\restrictedto\gamma$ and, for $x\in
D^p_\gamma$, $s^{p'}_x=s^p_x$, $F^{p'}_x=F^p_x$ and
$w^{p'}_x=\max\{w^p_x,2\cdot\size{F^p_x}\}$. Let
$N=2\cdot\tsum\{w^{p'}_x\st x\in D^p_\gamma\}$. Applying
Lemma~\ref{lem:repair} to $p$, $p'$ and $N$, we get a condition
$q\leq p$ as in the lemma. By induction hypothesis, there is a
condition $q^*\in W_{<\gamma}$, $q^*\leq q\restrictedto \gamma$.
Then $q^*\cup q\restrictedto\{\gamma\}$ extends $q$ (by
Lemma~\ref{lem:easycompleteembedding}) and belongs to
$W_{<\gamma+1}$.
\end{proof}

\begin{lemma}\label{lem:ccc}
$\bN_Q$ satisfies ccc.
\end{lemma}

\begin{proof}
Let $W$ be the dense set given by Lemma~\ref{lem:preccc}. If
$A\subseteq W$ is uncountable, then thin $A$ out to an uncountable
set $A'\subseteq A$ such that
\begin{enumerate}
\item[(1)]
$\{\bar D^p\st p\in A'\}$ is a $\Delta$-system with root $u$;
\item[(2)]
For $\xi\in u$, there is an $l_\xi$ such that $l^p_\xi=l_\xi$ for
all $p\in A'$;
\item[(3)]
$\{D^p\st p\in A'\}$ is a $\Delta$-system with root $U$;
\item[(4)]
For $x\in U$, there are $s_x$ and $w_x$ such that $s^p_x=s_x$ and
$w^p_x=w_x$ for all $p\in A'$;
\item[(5)]
For each $U'\subseteq U$, there is a number $k_{U'}$ such that for
each $p\in A'$,
\[
\tsum\{\size{F^p_z}\st \text{ for some }x\in U',\,z\in D^p_{\leq
x}\}
    =k_{U'}.
\]
Note that, because $p\in W$, we have $2k_{U'}\leq
\tsum\{w^p_z\st\text{ for some }x\in U',\,z\in D^p_{\leq x}\}$.
\end{enumerate}
Let $p$ and $q$ be any two conditions in $A'$. Let
$\xi_0<\xi_1<\dots<\xi_{k-1}$ be the increasing enumeration of
$\bar D^p\cup\bar D^q$. We will inductively define conditions
$r_i\in{\bN}_{<\xi_i+1}$, $i<k$, so that
\begin{enumerate}
\item
$r_i$ is a common extension of $p\restrictedto(\xi_i+1)$ and
$q\restrictedto(\xi_i+1)$;
\item
For each $i<k-1$, $r_{i+1}\restrictedto \xi_{i+1}\leq r_i$.
\end{enumerate}
Set $r_{-1}=\emptyset$. When $\xi_i\not\in u$, then only one of
$\bar D^p$, $\bar D^q$ contains $\xi_i$. If $\xi_i\in \bar D^p\ssm
\bar D^q$, then let $r_i=r_{i-1}\cup p\restrictedto\{\xi_i\}$.
Then $r_i$ inherits from $r_{i-1}$ and $p\restrictedto\{\xi_i\}$
the properties needed for being a condition. It extends
$p\restrictedto(\xi_i+1)$ by
Lemma~\ref{lem:easycompleteembedding}. It extends
$q\restrictedto(\xi_i+1)$ because the inclusion of the domains
holds and $q\restrictedto(\xi_i+1)=q\restrictedto(\xi_{i-1}+1)$,
so the relevant values of $x$ and $\xi$ for which $x\in D^q$ or
$\xi\in \bar D^q$ in clauses 6--9 of
Definition~\ref{def:iteration} applied to $r_i\leq
q\restrictedto(\xi_i+1)$ all have rank at most $\xi_{i-1}$ and
hence the clauses hold because $r_{i-1}\leq
q\restrictedto(\xi_{i-1}+1)$. Similarly if $\xi_i\in \bar D^q\ssm
\bar D^p$.

Now suppose $\xi_i=\gamma\in u$. Proceed as follows.
\begin{enumerate}
\item[(a)]
Let $L=\tsum\{w^p_x\st x\in D^p_\gamma\}
    +\tsum\{w^q_x\st x\in D^q_\gamma\}+l_\gamma$.
\item[(b)]
Get $r^*\in {\bN}_{<\gamma_i}$, $r^*\leq r_{i-1}$ which decides
the values of $\dot f\restrictedto L$ for $f\in F^p_x$, $x\in
D^p_\gamma$ and $f\in F^q_x$, $x\in D^q_\gamma$. For $n\in L$, let
$K_{x,n}$ be the set such that
\begin{enumerate}
\item[(i)]
$r^*\Vdash \{\dot f(n)\st \dot f\in F^p_x\}=K_{x,n}$, if $x\in
D^p_\gamma\ssm D^q_\gamma$;
\item[(ii)]
$r^*\Vdash \{\dot f(n)\st \dot f\in F^q_x\}=K_{x,n}$, if $x\in
D^q_\gamma\ssm D^p_\gamma$;
\item[(iii)]
$r^*\Vdash \{\dot f(n)\st \dot f\in F^p_x\cup F^q_x\}=K_{x,n}$, if
$x\in D^p_\gamma\cap D^q_\gamma$.
\end{enumerate}
Note that
\begin{align*}
\size{\tbigcup\{K_{x,n}\st x\in D^p_\gamma\cup D^q_\gamma\}}
 & \leq  \tsum\{w^p_x\st x\in D^p_\gamma\}
    + \tsum\{w^q_x\st x\in D^q_\gamma\}\\
 & \leq  2\cdot\max(\tsum\{w^p_x\st x\in D^p_\gamma\},
    \tsum\{w^q_x\st x\in D^q_\gamma\})\\
 & \leq  l_\gamma
\end{align*}
where the last inequality holds because $p,q\in W$.
\item[(c)]
For $n$ such that $l_\gamma\leq n<L$, define $s_x(n)$ as follows.
\begin{enumerate}
\item[(i)]
$s_x(n)
    =\tbigcup\{K_{z,n}\st z\in D^p_{\leq x}
        \text{ or for some } z'\in D^p_\gamma\cap D^q_\gamma,
        \,z\in (D^p\cup D^q)_\gamma
        \text{ and } z\leq z'\leq x\}$,
    if $x\in D^p_\gamma\ssm D^q_\gamma$;
\item[(ii)]
$s_x(n)
    =\tbigcup\{K_{z,n}\st z\in D^q_{\leq x}
        \text{ or for some } z'\in D^p_\gamma\cap D^q_\gamma,
        \,z\in (D^p\cup D^q)_\gamma\text{ and } z\leq z'\leq x\}$,
    if $x\in D^q_\gamma\ssm D^p_\gamma$;
\item[(iii)]
$s_x(n)=\tbigcup\{K_{z,n}\st z\in (D^p\cup D^q)_{\leq x}\}$, if
$x\in D^p_\gamma\cap D^q_\gamma$.
\end{enumerate}
Suppose $E\subseteq D^p_\gamma$ is downward closed. Then
\begin{align*}
\tbigcup\{s_x(n)\st x\in E\}
    = {} & \tbigcup\{K_{z,n}\st z\in(D^p\cup D^q)_{\leq x}
        \text{ for some } x\in E\cap U\}\\
    & \cup \tbigcup\{K_{z,n}\st z\in E
        \text{ and for no } x\in E\cap U\text{ do we have } z\leq x\}.
\end{align*}
So
\begin{align*}
\size{\tbigcup\{s_x(n)\st x\in E\}} \leq {}&
\tsum\{\size{F^p_z}\st z\in D^p_{\leq x}
       \text{ for some } x\in E\cap U\}\\
       & + \tsum\{\size{F^q_z}\st z\in D^q_{\leq x}
       \text{ for some } x\in E\cap U\}\\
       & + \tsum\{\size{F^p_z}\st z\in E
       \text{ and for no }x\in E\cap U\text{ do we have } z\leq x\}\\
\leq {}& 2k_{E\cap U} \\
       & + \tsum\{\size{F^p_z}\st z\in E
       \text{ and for no } x\in E\cap U\text{ do we have } z\leq x\}\\
\leq {}& \tsum\{w^p_z\st z\in D^p_{\leq x}
       \text{ for some } x\in E\cap U\}\\
       & + \tsum\{w^p_z\st z\in E
       \text{ and for no } x\in E\cap U\text{ do we have } z\leq x\}\\
   = {}& \tsum\{w^p_z\st z\in E\}.
\end{align*}
Similarly, if $E$ is a downward closed subset of $D^q_\gamma$,
then $\size{\tbigcup\{s_x(n)\st x\in E\}}\leq \tsum\{w^q_z\st z\in
E\}$.
\item[(d)] Let
$r_i=r^*\cup\{(s_x,w_x,F_x)\st x\in D^p_\gamma\cup D^q_\gamma\}$,
where the triples $(s_x,w_x,F_x)$ are obtained as follows.
\begin{enumerate}
\item[(i)]
Each $s_x$ has domain $L$, $s_x\restrictedto l_\gamma=s^p_x$ if
$x\in D^p_\gamma$ and $s_x\restrictedto l_\gamma=s^q_x$ if $x\in
D^q_\gamma$. (This is unambiguous if both clauses hold because of
item (4) in the list of properties of $A'$.) For $l_\gamma\leq
n<L$, $s_x(n)$ is as defined in (c).
\item[(ii)]
We have $w_x=w^p_x$ if $x\in D^p_\gamma$ and $w_x=w^q_x$ if $x\in
D^q_\gamma$ (and this is unambiguous if both clauses hold).
\item[(iii)]
For $x\in D^p\ssm D^q$, $F_x=F^p_x$. For $x\in D^q\ssm D^p$,
$F_x=F^q_x$. For $x\in D^p\cap D^q$, $F_x=F^p_x\cup F^q_x$.
\end{enumerate}
\end{enumerate}
We must check that $r_i$ is as desired. First we check that $r_i$
is a well-defined condition. In Definition 4.1, clause 1 and the
first and third statements of clause 2 hold by definition. The
second statement holds below rank $\xi_i$ because $r^*$ is a
condition. At rank $\gamma=\xi_i$, it holds because for each $x\in
(D^p\cup D^q)_\gamma$ and $n<L$, if $n<l_\gamma$ then
$\size{s_x(n)}\leq n$ because $p$ and $q$ are conditions and if
$l_\gamma\leq n<L$ then the argument at the end of (b) above shows
that $\size{s_x(n)}\leq l_\gamma\leq n$. For the last statement,
we have that $\size{F_x}$ is bounded by one of $\size{F^p_x}$,
$\size{F^q_x}$, $\size{F^p_x}+\size{F^q_x}$. In all cases, because
$p,q\in W$, we have that $\size{F_x}$ is bounded by either
$2\cdot\size{F^p_x}\leq w^p_x=w_x$ or $2\cdot\size{F^q_x}\leq
w^q_x=w_x$. For clause 3, the property is inherited from $r^*$ if
the rank of $x$ is less than $\xi_i$, and, if the rank of $x$ is
$\xi_i$, is inherited from $p$ or $q$ if $\xi_i\in\bar D^p\ssm
\bar D^q$ or $\xi_i\in\bar D^q\ssm \bar D^p$. Otherwise we have
$\tsum\{w_x\st x\in(D^p\cup D^q)_{\leq x}\}
    \leq \tsum\{w^p_x\st x\in D^p_{\leq x}\}
    + \tsum\{w^q_x\st x\in D^q_{\leq x}\}\leq l_{\xi_i}\leq L$.
Clause 4 is inherited from $r^*$ at ranks below $\xi_i$ and holds
by definition at rank $\xi_i$.

Now we check that $r$ extends $p$ and $q$. By symmetry, it its
enough to check that $r$ extends $p$. All of the clauses 5--9 in
the definition hold below rank $\xi_i$ because $r^*\leq
r_{i-1}\leq p\restrictedto\xi_{i-1}+1$. Consider now what they say
at rank $\gamma=\xi_i$. The inclusion of the domains and all but
the last part of 6 hold by definition of $r$. The last part of 6
holds because if $x\in D^p_\gamma$, $\dot f\in F^p_x$ and
$l_\gamma\leq n<L$, we chose $r^*$ so that
$r^*\Vdash_{{\bN}_{<\gamma}}\dot f(n)\in K_{x,n}\subseteq s_x(n)$.
Because $\dot f$ is a ${\bN}_x$-name and ${\bN}_x$ is completely
embedded in ${\bN}_{<\gamma}$, it follows that $r^*\restrictedto
x=r_i\restrictedto x$ also forces $\dot f(n)\in s_x(n)$.

The proof of clause 7 is a case by case analysis. Suppose $x,y\in
D^p_\gamma$, $x<y$ and $l_\gamma\leq n<L$. Each of $x$ and $y$
comes under either (c)(i) or (c)(iii). Since the formulas used
there are increasing functions of $x$, we need only consider the
following two cases.

Case 1. $x\in D^p_\gamma\ssm D^q_\gamma$ and $y\in D^p_\gamma\cap
D^q_\gamma$. Let $m\in s_x(n)$ and fix $z$ witnessing this. (So,
in particular, $m\in K_{z,n}$.)  We will show that
$K_{z,n}\subseteq s_y(n)$. If $z\in D^p_{\leq x}$, then also $z\in
D^p_{\leq y}$, so $K_{z,n}\subseteq s_y(n)$. The other possibility
is that for some $z'\in D^p_\gamma\cap D^q_\gamma$, $z\in (D^p\cup
D^q)_\gamma$ and $z\leq z'\leq x$. Then $z'\in (D^p\cup D^q)_{\leq
y}$, so again $K_{z,n}\subseteq s_y(n)$.

Case 2. $x\in D^p_\gamma\cap D^q_\gamma$ and $y\in D^p_\gamma\ssm
D^q_\gamma$. Fix $z\in (D^p\cup D^q)_{\leq x}$. Taking $z'=x$, we
have $z\leq z'<y$ witnessing that $K_{z,n}\subseteq s_y(n)$.

For clause 8, we have that $\tsum\{w_x\st x\in (D^p\cup
D^q)_\gamma\}
    \leq\tsum\{w^p_\gamma\st x\in D^p_\gamma\}
    +\tsum\{w^q_\gamma\st x\in D^q_\gamma\}=L-l_\gamma$
by the definition of $L$ in (a). Finally, clause 9 was checked in
(c).

For $i=k-1$, we get that $r_i$ is a common extension of $p$ and
$q$.

This completes the proof that ${\bN}_Q$ is ccc.
\end{proof}

\section{Proof of the main theorem}\label{sec:main}

This section is devoted to the proof of Hechler's theorem for the
null ideal. We will show that the forcing notion $\bN_Q$ satisfies
all the requirements of the theorem.

\begin{lemma}\label{lem:prolong}
For a downward closed set $B\subseteq Q$, $p\in\bN_Q$,
$\xi\in\bar{D^p}$ and $N<\omega$, there is $q\in\bN_B$ such that
$q\leq p$ and $l^p_\xi\geq N$.
\end{lemma}

\begin{proof}
Just apply Lemma~\ref{lem:repair} to $p'=p$ and $N$.
\end{proof}

\begin{lemma}\label{lem:join}
For a downward closed set $B\subseteq Q$, $p\in\bN_B$ and $a\in
B$, there is $q\in\bN_B$ such that $q\leq p$ and $a\in D^q$.
\end{lemma}

\begin{proof}
We may assume that $a\notin D^p$. Let $\alpha=\rk{a}$.

If $\alpha\notin\bar{D^p}$, then define $q\in\bN_B$ by letting
$D^q=D^p\cup\{a\}$, $s^q_a=\emptyset$, $w^q_a=0$,
$F^q_a=\emptyset$ and other components of $q$ are the same as $p$.

Now we assume that $\alpha\in\bar{D^p}$. Define an
$\alpha$-preextension $p'$ of $p$ in $\bN_B$ by letting
$D^{p'}=D^p\cup\{a\}$, $s^{p'}_a$ is arbitrary with length
$l^p_\alpha$, $w^{p'}_a=0$, $F^{p'}_a=\emptyset$ and other
components of $p'$ are the same as $p$. Apply
Lemma~\ref{lem:repair} to $p$, $p'$ and $N=0$, and we get
$q\in\bN_B$ with $q\leq p$ and $a\in D^q$.
\end{proof}

\begin{lemma}\label{lem:preaddname}
For a downward closed set $B\subseteq Q$, $p\in\bN_B$ and $a\in
D^p$, there is $q\in\bN_B$ such that $q\leq p$ and
$w^q_a\geq\size{F^q_a}+1$.
\end{lemma}

\begin{proof}
Let $\alpha=\rk{a}$. Define an $\alpha$-preextension $p'$ of $p$
in $\bN_B$ by letting $D^{p'}=D^p$, $w^{p'}_a=w^p_a+1$ and other
components of $p'$ are the same as $p$. Apply
Lemma~\ref{lem:repair} to $p$, $p'$ and $N=0$, and we get
$q\in\bN_B$ as required.
\end{proof}

\begin{lemma}\label{lem:addname}
For a downward closed set $B\subseteq Q$, $p\in\bN_B$, $a\in D^p$
and an $\bN_a$-name $\dot{f}$ for a function in $\omega^\omega$,
there is $q\in\bN_B$ such that $q\leq p$ and $\dot{f}\in F^q_a$.
\end{lemma}

\begin{proof}
First use Lemma~\ref{lem:preaddname}, and then put $\dot{f}$ into
$F^q_a$.
\end{proof}

Let $\V$ be a ground model and $G$ an $\bN_Q$-generic filter over
$\V$. For $a\in Q$, let $G\restrictedto
a=G\cap\bN_a=\{p\restrictedto a\st p\in G\}$. Then $G\restrictedto
a$ is an $\bN_a$-generic filter over $\V$.

In $\V[G]$, for $a\in Q$ let $\varphi_a=\tbigcup\{s^p_a\st p\in
G\text{ and }a\in D^p\}$. By Lemmata~\ref{lem:prolong} and
\ref{lem:join}, $\varphi_a$ is defined for every $a\in Q$, and
belongs to $\cS$.

\begin{lemma}\label{lem:localize}
In $\V[G]$, for every $a\in Q$ and
$f\in\omega^\omega\cap\V[G\restrictedto a]$, for all but finitely
many $n<\omega$ we have $f(n)\in\varphi_a(n)$.
\end{lemma}

\begin{proof}
Follows from Lemma~\ref{lem:addname} and the definition of
$\bN_Q$.
\end{proof}

\begin{lemma}\label{lem:cannibal}
For $a,b\in Q$, if $a<b$ and $\rk{a}=\rk{b}$, then for all but
finitely many $n<\omega$ we have
$\varphi_a(n)\subseteq\varphi_b(n)$.
\end{lemma}

\begin{proof}
Clear from the definition of $\bN_Q$.
\end{proof}

For $a\in Q$, let $H_a=H_{\varphi_a}$. Then each $H_a$ is a null
subset of $2^\omega$. We will show that, in $\V[G]$, the set
$\{H_a\st a\in Q\}$ is order-isomorphic to $(Q,\leq)$ and cofinal
in $(\N,\subseteq)$.

\begin{lemma}\label{lem:amoeba}
Let $a\in Q$. For a Borel null set $X\subseteq 2^\omega$ which is
coded in $\V[G\restrictedto a]$, we have $X\subseteq H_a$.
\end{lemma}

\begin{proof}
Follows from Lemma~\ref{lem:localize} and the observation in
Section~\ref{sec:combnull}.
\end{proof}

\begin{lemma}\label{lem:cofinal}
In $\V[G]$, for every null set $X\subseteq 2^\omega$ there is
$a\in Q$ satisfying $X\subseteq H_a$.
\end{lemma}

\begin{proof}
We may assume that $X$ is a Borel set in $\V[G]$. By our
assumption that countable subsets of $Q$ have strict upper bounds,
and because $\bN_Q$ is ccc, $X$ is coded in $\V[G\restrictedto a]$
for some $a\in Q$, and by Lemma~\ref{lem:amoeba}, we have
$X\subseteq H_a$.
\end{proof}

\begin{lemma}\label{lem:orderpreserving}
For $a,b\in Q$, if $a\leq b$ then $H_a\subseteq H_b$.
\end{lemma}

\begin{proof}
If $a\ll b$, then $H_a$ is coded in $\V[G\restrictedto b]$ and
hence $H_a\subseteq H_b$ follows from Lemma~\ref{lem:amoeba}. If
$a<b$ and $\rk{a}=\rk{b}$, then it follows from
Lemma~\ref{lem:cannibal} and the observation in
Section~\ref{sec:combnull}.
\end{proof}

For each $a\in Q$, let $r_a=r_{\varphi_a}$ and $R_a=R_{\varphi_a}$
as defined in Section~\ref{sec:combnull}. As we observed in
Section~\ref{sec:loc}, we define an $\bN_Q$-name $\dot{r}_a$ for
$r_a$ so that, for $p\in\bN_Q$ if $a\in D^p$ and $\lh{s^p_a}=n$
then $p$ decides the value of $\dot{r}_a\restrictedto n$.

\begin{lemma}\label{lem:incomp}
For $a,b\in Q$, if $a\not\leq b$ then $H_a\not\subseteq H_b$.
\end{lemma}

\begin{proof}
Suppose that $a\not\leq b$. Since we always have $R_b\cap
H_b=\emptyset$ and $R_b\neq\emptyset$, it suffices to show that
$R_b\subseteq H_a$.

Fix $p\in\bN_Q$ and $M<\omega$. By Lemmata~\ref{lem:join} and
\ref{lem:preaddname}, we may assume that $a,b\in D^p$ and
$w^p_a\geq\size{F^p_a}+1$.

We will find $q\leq p$ and $m>M$ which satisfy
$q\forces{\dot{r}_b(m)\in s^q_a(m)}$. This implies that for
infinitely many $m<\omega$ we have $r_b(m)\in\varphi_a(m)$, and
hence $R_b\subseteq H_a$.

Let $\alpha=\rk{a}$, $\beta=\rk{b}$, $B=\{x\in Q\st x\leq b\}$.
Note that $a\notin B$ by the assumption. Extend $p$ if necessary
to arrange the following.
\begin{center}
If $B_\alpha\not=\emptyset$, then $B_\alpha\cap
D^p\not=\emptyset$.
\end{center}
(The following observation is not used in the proof, but note for
clarity that because of the definition of rank for elements of
$Q\ssm R$, the ranks of the elements of a downward closed set need
not be an initial segment of the ordinals. For example, if
$R=\omega_1$ ordered as usual and $Q$ is $R$ with new elements
$e_\alpha$, where $e_\alpha\leq \alpha$ but no other relations
hold other than the ones needed to ensure transitivity, then
$e_\alpha$ has rank $\alpha$ and every subset of $\{e_\alpha\st
\alpha<\omega_1\}$ is downward closed. Thus the assumption
$B_\alpha\not=\emptyset$ can fail even if $\alpha<\beta$.)

We set $m=\max\{M,l^p_\alpha\}+1$.

Using Lemma~\ref{lem:prolong}, get $p^*\in{\bN}_B$ extending
$p\restrictedto B$ such that $\lh{s^{p^*}_b}\geq m+1$. By the
choice of $\dot{r}_b$, $p^*$ decides the value of $\dot{r}_b(m)$,
so let $k$ be such that $p^*\Vdash_{{\bN}_B}\dot{r}_b(m)=k$.

We will construct $q\in\bN_Q$ satisfying $q\leq p$ and $q\leq
p^*$, using an argument similar to, but somewhat more difficult
than, the proof of Lemma~\ref{lem:completeembedding}.

The proof which follows is really two similar but different
proofs, one for the case where $B_\alpha\not=\emptyset$ and one
for the case $B_\alpha=\emptyset$. In order to be able to write as
much as possible of the two proofs as one, we will use the abuse
of notation $\max\{l^{p^*}_\alpha,l^p_\alpha\}$ to designate
$l^{p^*}_\alpha$ when $B_\alpha\not=\emptyset$ and $l^p_\alpha$
when $B_\alpha=\emptyset$ (in which case $l^{p^*}_\alpha$ is
actually not defined).

We will be done if we build $q\leq p$ with $k\in s^q_a(m)$. For
$x\in D^{p^*}_\alpha$, let
$(s_x,w_x,F_x)=(s^{p^*}_x,w^{p^*}_x,F^{p^*}_x)$. For $x\in
D^p_\alpha\ssm D^{p^*}_\alpha$, let
$(s_x,w_x,F_x)=(s^{p}_x,w^{p}_x,F^{p}_x)$. Let
\[
L=\tsum\{w_x\st  x\in D^p_\alpha\cup D^{p^*}_\alpha\}
    +\max\{l^{p^*}_\alpha,l^p_\alpha\} + m + 1.
\]
Choose $q_0\in{\bN}_\alpha$ so that $q_0\leq
p\restrictedto\alpha$, $q_0\leq p^*\restrictedto\alpha$ (and hence
also $q_0\restrictedto B_{<\alpha}\leq p^*\restrictedto\alpha$),
and $q_0$ decides the values of $\dot{f}\restrictedto L$ for all
$\dot{f}\in\tbigcup\{F_x\st  x\in D^p_\alpha\cup
D^{p^*}_\alpha\}$. For $x\in D^p_\alpha\cup D^{p^*}_\alpha$ and
$n\in L\ssm\lh{s_x}$, let $K_{x,n}\subseteq\omega$ be the set
satisfying $q_0\Vdash{K_{x,n}=\{\dot{f}(n)\st \dot{f}\in F_x\}}$.
For $x\in D^p_\alpha\cup D^{p^*}_\alpha$ and $n\in L\ssm\lh{s_x}$,
if $(x,n)\neq (a,m)$ then let $K'_{x,n}=K_{x,n}$, and let
$K'_{a,m}=K_{a,m}\cup\{k\}$. By the assumption that
$w^p_a\geq\size{F^p_a}+1$, we have $\size{K'_{x,n}}\leq w_x$ for
all $x\in D^p_\alpha\cup D^{p^*}_\alpha$ and $n\in L\ssm\lh{s_x}$.

Define $s^*_x$ for $x\in D^p_\alpha\cup D^{p^*}_\alpha$ as
follows. If $x\in D^{p^*}_\alpha$, then $\lh{s^*_x}=L$,
$s^*_x\restrictedto l^{p^*}_\alpha=s_x$, and for $n\in L\ssm
l^{p^*}_\alpha$,
\[
s^*_x(n)=\tbigcup\{K'_{z,n}\st  z\in D^{p^*}_{\leq x}\}.
\]
If $x\in D^p_\alpha\ssm D^{p^*}_\alpha$, then $\lh{s^*_x}=L$,
$s^*_x\restrictedto l^p_\alpha=s_x$, and for $n\in L\ssm
l^p_\alpha$,
\[
s^*_x(n)=
\begin{cases}
\tbigcup\{s_z(n)\st  z\in D^p_{\leq x}\cap D^{p^*}_\alpha\}
   \cup\tbigcup\{K'_{z,n}\st  z\in D^p_{\leq x}\ssm
   D^{p^*}_\alpha\},
& l^p_\alpha\leq
   n<\max\{l^{p^*}_\alpha,l^p_\alpha\}\\
\tbigcup\{K'_{z,n}\st z\in (D^p_\alpha\cup D^{p^*}_\alpha)_{\leq
x}\},
   & \max\{l^{p^*}_\alpha,l^{p}_\alpha\}\leq n<L
\end{cases}
\]
Define $q_1$ by $q_1=\{(s^{q_1}_x,w^{q_1}_x,F^{q_1}_x)
    \st  x\in D^{q_0}\cup D^{p^*}\cup D^p_\alpha\}$
where
\begin{enumerate}
\item For $x\in D^{q_0}$,
$(s^{q_1}_x,w^{q_1}_x,F^{q_1}_x)=(s^{q_0}_x,w^{q_0}_x,F^{q_0}_x)$
\item For $x\in D^p_\alpha\cup D^{p^*}_\alpha$,
$(s^{q_1}_x,w^{q_1}_x,F^{q_1}_x)=(s^{*}_x,w_x,F_x)$
\item For $x\in D^{p^*}\ssm Q_{<\alpha+1}$,
$(s^{q_1}_x,w^{q_1}_x,F^{q_1}_x)=(s^{p^*}_x,w^{p^*}_x,F^{p^*}_x)$
\end{enumerate}
We now check that $q_1\in{\bN}_Q$. The requirements of
Definition~\ref{def:iteration} are satisfied below (resp.\ above)
rank $\alpha$ because $q_0$ (resp.\ $p^*$) is a condition.
Consider what they say at rank $\alpha$. The first clause is
trivial. The fourth holds because the $s^{q_1}_x$'s all have
domain $L$. The third clause can be checked in two cases.
\begin{enumerate}
\item[(i)]
If $x\in D^{p^*}_\alpha$, then $D^{q_1}_{\leq x}=(D^p\cup
D^{p^*})_{\leq x}=D^{p^*}_{\leq x}$, so $\tsum\{w^{q_1}_z\st z\in
D^{q_1}_{\leq x}\}
    = \tsum\{w^{p^*}_z\st z\in D^{p^*}_{\leq x}\}
    \leq l^{p^*}_\alpha\leq L$.
\item[(ii)]
If $x\in D^{p}_\alpha\ssm D^{p^*}_\alpha$, then $D^{q_1}_{\leq
x}=D^{p}_{\leq x}\cup D^{p^*}_{\leq x}$, so $ \tsum\{w^{q_1}_z\st
z\in D^{q_1}_{\leq x}\}
    = \tsum\{w_z\st z\in D^{p}_{\leq x}\cup D^{p^*}_{\leq x}\}
    \leq \tsum\{w_z\st z\in D^{p}_{\alpha}\cup D^{p^*}_{\alpha}\}
    \leq L.$
\end{enumerate}
For the second, all the requirements except that the $s^{q_1}_x$'s
are partial slaloms follow from the fact that $p$ and $p^*$ are
conditions. We need to check that $\size{s^*_x(n)}\leq n$ for each
relevant $n$. If $x\in D^{p^*}_\alpha$, then for
$l^{p^*}_\alpha\leq n<L$, we have
$\size{s^*_x(n)}\leq\tsum\{w^{p^*}_z\st z\in D^{p^*}_{\leq x}\}
    \leq\lh{s^{p^*}_x}=l^{p^*}_\alpha\leq n$.
If $x\in D^p_\alpha\ssm D^{p^*}_\alpha$, we consider three cases.

Case 1. $l^p_\alpha\leq n<\max\{l^{p^*}_\alpha,l^p_\alpha\}$. In
order for this case to be non-vacuous, we must have $\alpha\in\bar
D^{p\restrictedto B}$. Then Definition~\ref{def:iteration}(9) for
$p^*\leq p\restrictedto B$ with $E=D^p_{\leq x}\cap D^{p^*}$ gives
\begin{align*}
\size{s^*_x(n)}
    & \leq  \tsum\{w^{p}_z\st z\in E\} + (n-l^p_\alpha) +
        \tsum\{w^{p}_z\st z\in D^{p}_{\leq x}\ssm E\}\\
    & =   \tsum\{w^{p}_z\st z\in D^{p}_{\leq x}\} +
        (n-l^p_\alpha)\\
    & \leq  l^p_\alpha + (n - l^p_\alpha) = n.
\end{align*}

Case 2. $\max\{l^{p^*}_\alpha,l^p_\alpha\}\leq n<L$. If $\alpha\in
\bar{D}^{p\restrictedto B}$, then
Definition~\ref{def:iteration}(8) for $p^*\leq p\restrictedto B$
gives
\[
\tsum\{w^{p^*}_z\st z\in D^{p^*}_\alpha\}
    \leq \tsum\{w^p_z\st z\in D^{p\restrictedto B}_\alpha\}
    +(l^{p^*}_\alpha-l^p_\alpha).
\]
Removing terms with $z\not\leq x$ from both sides (see
Remark~\ref{rem:rwgrowthdiscardterms}) gives
\[
\tsum\{w^{p^*}_z\st z\in D^{p^*}_{\leq x}\}
    \leq \tsum\{w^p_z\st z\in D^p_{\leq x}\cap B\}
        +(l^{p^*}_\alpha-l^p_\alpha) .
\]
From the formula for $s_x(n)$ we now get
\begin{align*}
\size{s_x(n)} & \leq  \tsum\{w^{p^*}_z\st z\in D^{p^*}_{\leq x}\}
    + \tsum\{w^{p}_z\st z\in D^p_{\leq x}\ssm B\}\\
& \leq  \tsum\{w^p_z\st z\in D^p_{\leq x}\cap B\}
    + (l^{p^*}_\alpha-l^p_\alpha)
    + \tsum\{w^{p}_z\st z\in D^p_{\leq x}\ssm B\}\\
& =  \tsum\{w^p_z\st z\in D^p_{\leq x}\}
    + (l^{p^*}_\alpha-l^p_\alpha)\\
& \leq  l^p_\alpha + (l^{p^*}_\alpha - l^p_\alpha)
    = l^{p^*}_\alpha\leq n.
\end{align*}
If $\alpha\not\in \bar D^{p\restrictedto B}$, then
$B_\alpha=\emptyset$, so $\alpha\not\in \bar D^{p^*}$. The formula
for $s^*_x(n)$ thus reduces to $s^*_x(n)=\tbigcup\{K'_{z,n}\st
z\in D^p_{\leq x}\}$, and hence $\size{s^*_x(n)}\leq
\tsum\{w^p_z\st z\in D^p_{\leq x}\}
    \leq l^p_\alpha\leq n$.

Thus, $q_1$ is a condition. We now check
Definition~\ref{def:iteration}(5--9) for $q_1\leq p^*$ and
$q_1\leq p\restrictedto B\cup Q_{<\alpha+1}$. (We only need the
latter, but the former is needed at one point of the proof.)
Clause~5 follows from the definition of $q_1$. For clauses~6--9,
first note that below rank $\alpha$, they hold because $q_0\leq
p\restrictedto\alpha$ and $q_0\leq p^*\restrictedto\alpha$.
Consider what happens at rank $\alpha$. Clause~6 holds because for
$x\in D^{p}_\alpha\cup D^{p^*}_\alpha$ and all the relevant values
of $\dot f$ and $n$, we have from the definitions that
$q_0\Vdash\dot f(n)\in K_{x,n}$ and $K_{x,n}\subseteq s^*_x(n)$.
For clause~7, we consider three cases. Let $x<y$ be elements of
$D^p_\alpha\cup D^{p^*}_\alpha$.
\begin{enumerate}
\item[(i)]
If $x,y\in D^{p^*}_\alpha$, then for checking $q_1\leq p^*$, just
use the monotonicity of $s^*_x(n)$ as a function of $x$. For
checking $q_1\leq p\restrictedto B\cup Q_{<\alpha+1}$, we also
need to consider values of $n$ such that $l^p_\alpha\leq
n<l^{p^*}_\alpha$. But then $s^*_x(n)=s^{p^*}_x(n)\subseteq
s^{p^*}_y(n) = s^*_y(n)$ because $p^*\leq p\restrictedto B$.

This is the only case to consider for checking clause~7 for
$q_1\leq p^*$ at stage $\alpha$. The remaining cases deal with
checking $q_1\leq p\restrictedto B\cup Q_{<\alpha+1}$. Note that
if $y\in D^{p^*}_\alpha\cap D^p_\alpha=D^{p}_\alpha\cap B$ then
also $x \in D^{p^*}_\alpha\cap D^p_\alpha$ since $B$ is downward
closed.
\item[(ii)]
If $x,y\in D^p_\alpha\ssm D^{p^*}_\alpha$, use the monotonicity of
$s^*_x(n)$ as a function of $x$.
\item[(iii)]
If $x\in D^{p^*}_\alpha\cap D^p_\alpha$ and $y\in D^p_\alpha\ssm
D^{p^*}_\alpha$, then consider first a value of $n$ such that
$l^p_\alpha\leq n<l^{p^*}_\alpha$. We have $s^*_x(n)=s_x(n)
    \subseteq \tbigcup\{s_z(n)\st z\in D^p_{\leq y}\cap D^{p^*}_\alpha\}
    \subseteq s^*_y(n)$.
Next consider $n$ such that $l^{p^*}_\alpha\leq n<L$. We have
$s^*_x(n)
    = \tbigcup\{K'_{z,n}\st z\in D^{p^*}_{\leq x}\}
    \subseteq \tbigcup\{K'_{z,n}\st
        z\in (D^p_\alpha\cup D^{p^*}_\alpha)_{\leq y}\}
    =s^*_y(n)$.
\end{enumerate}
That takes care of clause~7. Clause~8 follows from the fact that
if $\alpha\in \bar D^{p^*}$, then from the definition of $L$ we
have $\tsum\{w_x\st x\in D^{p^*}_\alpha\cup D^p_\alpha\}\leq
L-l^{p^*}_\alpha$, and if $\alpha\in \bar D^p_\alpha\ssm \bar
D^{p^*}_\alpha$, then $\tsum\{w_x\st x\in D^{p^*}_\alpha\cup
D^p_\alpha\}
    \leq L-l^{p}_\alpha$.
For clause 9, first we check $q_1\leq p^*$. If $\alpha\in \bar
D^{p^*}$, $E\subseteq D^{p^*}_\alpha$ is downward closed in
$D^{p^*}_\alpha$ and $l^{p^*}_\alpha\leq n<L$, then
$\size{\tbigcup\{s^*_x(n)\st x\in E\}}
    =\size{\tbigcup\{K'_{x,n}\st x\in E\}}
    \leq\tsum\{w^{p^*}_x\st x\in E\}$.
Next we check $q_1\leq p\restrictedto B\cup Q_{<\alpha+1}$. Note
that the elements of rank $\alpha$ are the same for the domains of
$p$ and $p\restrictedto B\cup Q_{<\alpha+1}$. Also $\alpha\in
D^p_\alpha$ since $a\in D^p$. Let $E\subseteq D^p_\alpha$ be
downward closed. Consider two cases.

Case 1. $l^p_\alpha\leq n<l^{p^*}_\alpha$. We have
\begin{align*}
\size{\tbigcup\{s^*_x(n)\st x\in E\}}
    &= \size{\tbigcup\{s^*_x(n)\st x\in E\cap B\}
        \cup\tbigcup\{s^*_x(n)\st x\in E\ssm B\}}   \\
    &= \size{\tbigcup\{s^{p^*}_x(n)\st x\in E\cap B\}
        \cup \tbigcup\{K'_{x,n}\st x\in E\ssm B\}}  \\
    &\leq \tsum\{w^{p\restrictedto B}_x\st x\in E\cap B\}+(n-l^p_\alpha)
        +\tsum\{w^p_x\st x\in E\ssm B\} \\
    &= \tsum\{w^{p}_x\st x\in E\}+(n-l^p_\alpha).
\end{align*}

Case 2. $\max\{l^{p^*}_\alpha,l^p_\alpha\}\leq n<L$. Let
$E'=\{z\in D^{p^*}_\alpha\st \text{ for some }x\in E\text{, }z\leq
x\}$. We have
\begin{align*}
\size{\tbigcup\{s^*_x(n)\st x\in E\}}
    &= \size{\tbigcup\{s^*_x(n)\st x\in E\cap B\}
        \cup\tbigcup\{s^*_x(n)\st x\in E\ssm B\}}   \\
    &= \size{\tbigcup\{K'_{x,n}\st x\in E'\}
        \cup \tbigcup\{K'_{x,n}\st x\in E\ssm B\}}  \\
    &\leq \tsum\{w^{p^*}_x\st x\in E'\} + \tsum\{w^p_x\st x\in E\ssm B\}.
\end{align*}
If $E'$ is empty, then this last expression is
$\leq\tsum\{w^p_x\st x\in E\}$. If not, then
Definition~\ref{def:iteration}(8) applied to $p^*\leq
p\restrictedto B$ (with terms outside $E'$ eliminated from both
sides) gives that
\begin{align*}
&\tsum\{w^{p^*}_x\st x\in E'\} + \tsum\{w^p_x\st x\in E\ssm B\} \\
    &\leq \tsum\{w^p_x\st x\in E\cap B\}
        + (l^{p^*}_\alpha - l^p_\alpha)
        + \tsum\{w^p_x\st x\in E\ssm B\}    \\
    &= \tsum\{w^p_x\st x\in E\} + (l^{p^*}_\alpha - l^p_\alpha) \\
    &\leq \tsum\{w^p_x\st x\in E\} + (n - l^p_\alpha).
\end{align*}

Thus, the conditions for $q_1\leq p^*$ and $q_1\leq p\restrictedto
B\cup Q_{<\alpha+1}$ hold up to rank $\alpha$. Above rank
$\alpha$, $q_1$ agrees with $p^*$, so
Definition~\ref{def:iteration}(6--9) hold trivially for $q_1\leq
p^*$. For $q_1\leq p\restrictedto B\cup Q_{<\alpha+1}$ we need to
prove the the clauses for $\xi>\alpha$. All of them follow from
the fact that $p^*\leq p\restrictedto B$, $q_1\restrictedto\xi\leq
p^*\restrictedto\xi$, and $q_1$ agrees with $p^*$ at rank $\xi$.
(The fact that $q_1\restrictedto\xi\leq p^*\restrictedto\xi$ is
used to check the last part of clause~6.)

Now we apply Lemma~\ref{lem:completeembedding} to $p$ and $q_1$,
and we get $q\in\bN_Q$ such that $q\leq p$ and
$q\forces{\dot{r}_b(m)\in s^q_a(m)}$.
\end{proof}

Now we have the following main theorem.

\begin{theorem}\label{thm:main}
Let $\N$ be the collection of null sets in $2^\omega$. Suppose
that $Q$ is a partially ordered set such that every countable
subset of $Q$ has a strict upper bound in $Q$. Then in any forcing
extension by $\bN_Q$, $(\N,\subseteq)$ contains a cofinal subset
$\{H_a\st a\in Q\}$ which is order-isomorphic to $(Q,\leq)$, that
is,
\begin{enumerate}
\item for every $X\in\N$ there is $a\in Q$
    such that $X\subseteq H_a$,
and
\item for $a,b\in Q$,
$H_a\subseteq H_b$ if and only if $a\leq b$.
\end{enumerate}
\end{theorem}

\subsection*{Acknowledgement}

The authors thank J.~Zapletal for pointing out the relevance of
the result of Hjorth mentioned in the introduction. They also
thank T.~Bartoszy\'nski, J.~Brendle, S.~Fuchino, S.~Kamo and
T.~Miyamoto for their helpful comments, suggestions and discussion
during this work.


\end{document}